\documentclass[11pt,dvips]{article}
\usepackage{amsmath,amsfonts,amssymb,amsthm,graphicx,verbatim,subfig}
\usepackage[all]{xy}
\usepackage[a4paper,left=20mm,right=20mm, top=15mm, bottom=15mm]{geometry}
\ifx\pdftexversion\undefined

\usepackage[a4paper,colorlinks,link
=black,filecolor=black,citecolor=black,urlcolor=black,pdfstartview=FitH]{hyperref}
\else

\usepackage[a4paper,colorlinks,linkcolor=black,filecolor=black,citecolor=black,urlcolor=black,pdfstartview=FitH]{hyperref}
\fi

\font\sixbb=msbm6
\font\eightbb=msbm8
\font\twelvebb=msbm10 scaled 1095
\newfam\bbfam
\textfont\bbfam=\twelvebb \scriptfont\bbfam=\eightbb
                           \scriptscriptfont\bbfam=\sixbb
\def\bb{\fam\bbfam\twelvebb}
\newcommand{\Rea}{{\bb R}}

\newcommand{\Rat}{{\bb Q}}

\newcommand{\FF}{{\bb F}}

\newtheorem{theorem}{\bf Theorem}[section]

\newtheorem{proposition}[theorem]{\bf Proposition}
\newtheorem{corollary}[theorem]{\bf Corollary}
\newtheorem{definition}[theorem]{\bf Definition}
\newcommand{\enp}{\begin{flushright} $\Box$ \end{flushright}}
\newcommand{\beq}[0]{\begin{equation}}
\newcommand{\enq}[0]{\end{equation}}

\newcommand{\thh}{\tilde{H}}

\newcommand{\cf}{{\cal F}}
\newcommand{\cg}{{\cal G}}
\newcommand{\ck}{\mathcal{K}}

\newcommand{\til}{\tilde{L}}
\newcommand{\lk}{{\rm lk}}

\newcommand{\sd}{{\rm sd}}

\newcommand{\str}{{\rm st}}

\title{Relative Leray Numbers via Spectral Sequences}
\begin{document}
\author{Gil Kalai \thanks{
Einstein Institute of Mathematics, Hebrew University, Jerusalem 91904, Israel, and
Efi Arazy School of Computer Science, IDC, Herzliya. e-mail: kalai@math.huji.ac.il~.
Supported by by  ISF grant 1612/17.}
\and Roy Meshulam\thanks{Department of Mathematics,
Technion, Haifa 32000, Israel. e-mail:
meshulam@technion.ac.il~. Supported by ISF grant 686/20.}}

\maketitle
\pagestyle{plain}

\begin{abstract}
Let $\FF$ be a fixed field and let $X$ be a simplicial complex on the vertex set $V$.
The Leray number $L(X;\FF)$ is the minimal $d$ such that for all $i \geq d$ and $S \subset V$, the induced complex $X[S]$ satisfies $\thh_i(X[S];\FF)=0$. Leray numbers play a role in formulating and proving topological Helly type theorems.
For two complexes $X,Y$ on the same vertex set $V$, define the relative Leray number $L_Y(X;\FF)$ as the minimal $d$ such that $\thh_i(X[V \setminus \sigma];\FF)=0$ for all $i \geq d$ and $\sigma \in Y$. In this paper we extend the topological colorful Helly theorem to the relative setting.
Our main tool is a spectral sequence for the intersection of complexes indexed by a geometric lattice.
\end{abstract}

\section{Introduction}
\label{s:intro}

Let $\FF$ be a fixed field and let $X$ be a simplicial complex on the vertex set $V$. All homology and cohomology appearing in the the sequel will be with $\FF$ coefficients. The induced subcomplex of $X$ on a subset $S \subset V$ is
$X[S]=\{\sigma \in X: \sigma \subset S\}$.
\begin{definition}
\label{d:leray}
The \emph{Leray number} $L(X)=L(X;\FF)$ of $X$ over $\FF$ is the minimal $d$ such that
$\thh_i(X[S])=0$
for all $S \subset V$ and $i \geq d$. The complex $X$ is \emph{$d$-Leray} over $\FF$ if $L(X) \leq d$.
\end{definition}
First introduced by Wegner \cite{Wegner75}, the family $\mathcal{L}^d=\mathcal{L}_{\FF}^d$ of $d$-Leray complexes
over the field $\FF$, has the following relevance to Helly type theorems.
Let $\cf$ be a family of sets. The \emph{Helly number}
$h(\cf)$ is the minimal positive integer $h$ such that if a
finite subfamily $\cg \subset \cf$ satisfies $\bigcap \cg' \neq
\emptyset$ for all $\cg' \subset \cg$ of cardinality $\leq h$, then
$\bigcap \cg \neq \emptyset$. Let $h(\cf)=\infty$ if no such finite
$h$ exists. For example, Helly's classical theorem asserts that the Helly number of the family of
convex sets in $\Rea^d$ is $d+1$.  Helly type theorems can often be formulated as properties of the associated nerves. Recall that the {\it nerve} of a family of sets $\cf$ is the simplicial complex $N(\cf)$ on the vertex set $\cf$, whose simplices are all subfamilies $\cg \subset \cf$ such that $\bigcap\cg \neq \emptyset$. A simple link between the Helly and Leray numbers is the inequality
$h(\cf) \leq L\left(N(\cf)\right)+1$ (see e.g. (1.2) in \cite{KM08}). A simplicial complex $X$ is {\it $d$-representable} if $X=N(\ck)$ for a family $\ck$ of convex sets in $\Rea^d$. Let $\mathcal{K}^d$ be the set of all $d$-representable complexes. Helly's theorem can then be stated as follows: If $X \in \mathcal{K}^d$ contain the full $d$-skeleton of its vertex set, then $X$ is a simplex.
The nerve lemma (see e.g. \cite{Bjorner03}) implies that
$\mathcal{K}^d \subset \mathcal{L}^d$, but the latter family is much richer, and there is substantial interest in understanding to what extent Helly type statements for $\mathcal{K}^d$ remain true for $\mathcal{L}^d$. A basic example is the following. A finite family $\cf$ of simplicial complexes in $\Rea^d$ is a {\it good cover} if for any $\cf' \subset \cf$, the intersection
$\bigcap \cf'$ is either empty or contractible. If $\cf$ is a good
cover in $\Rea^d$, then by the nerve lemma $N(\cf)$ is homotopic to $\bigcup \cf$ and therefore
$L\left(N(\cf)\right) \leq d$. Hence follows the Topological Helly's
Theorem: If $\cf$ is a good cover in $\Rea^d$, then $h(\cf) \leq L\left(N(\cf)\right)+1 \leq d+1$.

The Colorful Helly Theorem due to B\'{a}r\'{a}ny and Lov\'{a}sz  \cite{B82} is a fundamental result with a number of important applications in discrete geometry.
\begin{theorem}[\cite{B82}]
\label{t:chelly}
Let $\ck_1,\ldots,\ck_{d+1}$ be $d+1$ finite families of convex sets in $\Rea^d$, such that
$\bigcap_{i=1}^{d+1}K_i \neq \emptyset$ for all choices of $K_1 \in \ck_1, \ldots, K_{d+1} \in \ck_{d+1}$. Then there exists an $1 \leq i \leq d+1$ such that $\bigcap_{K \in \ck_i} K \neq \emptyset$.
\end{theorem}
\noindent
In \cite{KM05} we showed that the $d$-representability of $X=N(\bigcup_{i=1}^{d+1} \ck_i)$ can be replaced by the weaker assumption that $X$ is $d$-Leray.
\begin{theorem}[\cite{KM05}]
\label{t:tchelly}
Let $V=\bigcup_{i=1}^{d+1} V_i$ be a partition of $V$, and let $X$ be a $d$-Leray complex on $V$.
View each $V_i$ as a $0$-dimensional complex and suppose that $X$ contains the join $V_1* \cdots*V_{d+1}$.
Then there exists an $1 \leq i \leq d+1$ such that $V_i$ is a simplex of $X$.
\end{theorem}
\noindent
In fact, the transversal matroid
$V_1*\cdots *V_{d+1}$ in the statement of Theorem \ref{t:tchelly}, can be replaced by an arbitrary matroid. In the sequel we identify a matroid with the simplicial complex of its independent sets.
We recall that every induced subcomplex $M[S]$ of a matroid $M$ is pure, namely all maximal faces of $M[S]$ have the same dimension. This property can actually serve as the definition of a matroid in terms of the simplicial complex of its independent sets.
\begin{theorem}[\cite{KM05}]
\label{t:tmchelly}
Let $M$ be a matroid with a rank function $\rho_M$, and let $X$ be a $d$-Leray complex over some field $\FF$, both on the same vertex set $V$. If $M \subset X$, then there exists a $\sigma \in X$ such that $\rho_M(V \setminus \sigma) \leq d$.
\end{theorem}
In this paper we prove a generalization of Theorem \ref{t:tmchelly} using a new spectral sequence approach.
Let $X$ and $Y$ be two complexes on the same vertex set $V$.
\begin{definition}
\label{d:releray}
The \emph{relative Leray number} of $X$ with respect to $Y$ is
$$
L_Y(X)=L_Y(X;\FF)=\min\{d: \thh_{i}(X[V\setminus\sigma])=0~
\text{~for~all~} i \geq d \text{~and~} \sigma \in Y \}.$$
\end{definition}
\noindent
Our main result is the following relative extension of Theorem \ref{t:tmchelly}.
Let $X,Y$ be simplicial complexes on the vertex set $V$. Let $Y^{\vee}=\{A \subset V: V \setminus A \not\in Y\}$ denote the Alexander dual of $Y$ as a subcomplex of the simplex on $V$.
\begin{theorem}
\label{t:rtch}
Let $M$ be a matroid such that $Y^{\vee} \subset M \subset X$. Then there exists a simplex $\sigma \in X$ such that $\rho_M(V \setminus \sigma) \leq L_Y(X)$.
\end{theorem}

The paper is organized as follows. In section \ref{s:links} we  give a characterization of the relative Leray numbers in terms of links. In section \ref{s:interm} we construct a Mayer-Vietoris type spectral sequence (Proposition \ref{p:specs1}), and use it to establish a homological non-vanishing criterion (Corollary \ref{c:intz}) for certain families of complexes indexed by a geometric lattice. This result, which may be of independent interest, is the main ingredient in the proof of Theorem \ref{t:rtch} given in section \ref{s:rtch}.

\section{Relative Leray Numbers via Links}
\label{s:links}

We first recall a few definitions.
Let $X$ be a simplicial complex on the vertex set $V$.
The star and link of a subset $\tau \subset V$ are given by
\begin{equation*}
\begin{split}
\str(X,\tau)&=\{\sigma \in X: \sigma \cup \tau \in X\}, \\
\lk(X,\tau)&=\{\sigma \in \str(X,\tau): \sigma \cap \tau=\emptyset\}.
\end{split}
\end{equation*}
\noindent
Note that if $\tau \not\in X$, then $\str(X,\tau)=\lk(X,\tau)=\{\,\}$ is the void complex.
In particular, if $\thh_*(\lk(X,\tau)) \neq 0$ then $\tau \in X$.
It is well known that $L(X) \leq d$ iff $\thh_i\left(\lk(X,\sigma)\right)=0$ for all simplices $\sigma \in X$ and $i \geq d$.  The relative version of this fact is the following
\begin{proposition}
\label{p:lerayef}
Let $X,Y$ be complexes on the vertex set $V$. Then
\begin{equation*}
\label{e:lerayef}
L_Y(X)=\til_Y(X):=
\min\{d: \thh_i\left(\lk(X,\sigma)\right)=0 \text{~for~all~} i \geq d \text{~and~}  \sigma \in Y \}.
\end{equation*}
\end{proposition}
\noindent
Proposition \ref{p:lerayef} is implicit in the proof of Proposition 3.1 in \cite{KM06}, and
can also be deduced from a result of Bayer, Charalambous and Popescu (see Theorem 2.8 in \cite{BCP99}). For completeness, we include a simple direct proof of a slightly stronger result, following the argument in \cite{KM06}.
\begin{definition}
Let $X$ be a complex on the vertex set $V$ and let $A \subset V$.
The pair $(X,A)$ satisfies
{\em property $P_d(k_1,k_2)$}, if $\thh_i\bigl(\lk\left(X[V\setminus \sigma_1],\sigma_2\bigr) \right)=0$ for
all $i \geq d$ and all disjoint $\sigma_1,\sigma_2 \subset A$ such that $|\sigma_1| \leq k_1$, $|\sigma_2| \leq k_2$.
\end{definition}
\begin{proposition}
\label{p:ppk}
For a fixed pair $(X,A)$ and $k_1 \geq 0$, $k_2 \geq 1$, the properties $P_d(k_1,k_2)$ and $P_d(k_1+1,k_2-1)$ are equivalent.
\end{proposition}
\noindent
{\bf Proof.} Let $\tau_1,\tau_2$ be disjoint subsets of $A$ such that $|\tau_1| \leq k_1+1$, $|\tau_2| \leq k_2-1$ and suppose $v \in \tau_1$.
Let $\sigma_1=\tau_1 \setminus \{v\}$, $\sigma_2=\tau_2 \cup \{v\}$, and let
$$
~~~Z_1=\lk(X[V\setminus \tau_1],\tau_2)~~~~~,~~~Z_2=\str(\lk(X[V \setminus \sigma_1],\tau_2),v).
$$
\noindent
Then
$$
Z_1 \cup Z_2= \lk(X[V \setminus \sigma_1],\tau_2)~~~,~~~Z_1 \cap Z_2= \lk(X[V \setminus \sigma_1],\sigma_2).
$$
By Mayer-Vietoris there is an exact sequence
\begin{equation}
\label{e:mayerv}
\begin{split}
&\ldots  \rightarrow \thh_{i+1} \bigl(  \lk(X[V \setminus \sigma_1],\tau_2) \bigr)
\rightarrow \thh_i \bigl(  \lk(X[V \setminus \sigma_1],\sigma_2) \bigr) \\
&\rightarrow \thh_i \bigl( \lk(X[V\setminus \tau_1],\tau_2) \bigr) \rightarrow
\thh_i \bigl( \lk(X[V \setminus \sigma_1],\tau_2) \bigr) \rightarrow  \ldots ~~~.
\end{split}
\end{equation}
${\bf P_d(k_1,k_2) \Rightarrow P_d(k_1+1,k_2-1):}~~$ Suppose $(X,A)$ satisfies
$P_d(k_1,k_2)$ and let $i \geq d$. Let $\tau_1, \tau_2$ be disjoint subsets of $A$ such that $|\tau_1|=k_1+1$, $|\tau_2| \leq k_2-1$.
Choose $v \in \tau_1$, and let $\sigma_1=\tau_1\setminus \{v\}$, $\sigma_2=\tau_2 \cup\{v\}$.
The assumption that $P_d(k_1,k_2)$ holds implies that
the second and the fourth terms in (\ref{e:mayerv}) vanish. It
follows that $\thh_i \bigl( \lk(X[V\setminus \tau_1],\tau_2) \bigr)=0$
as required.
\\
${\bf P_d(k_1+1,k_2-1) \Rightarrow P_d(k_1,k_2):}~~$ Suppose $(X,A)$ satisfies
$P_d(k_1+1,k_2-1)$ and let $i \geq d$. Let $\sigma_1,\sigma_2$ be disjoint subsets of $A$ such that
$|\sigma_1| \leq k_1, |\sigma_2|=k_2$. Choose $v \in \sigma_2$, and let $\tau_1=\sigma_1 \cup\{v\}$ and
$\tau_2=\sigma_2 \setminus \{v\}$. The assumption that $P_d(k_1+1,k_2-1)$ holds implies that
the first and the third terms in (\ref{e:mayerv}) vanish. It
follows that $\thh_i \bigl( \lk(X[V\setminus \sigma_1],\sigma_2) \bigr)=0$
as required.
{\enp}
\noindent
{\bf Proof of Proposition \ref{p:lerayef}.} Clearly, $L_Y(X)=\max\{L_{A}(X): A \in Y\}$, and
$\til_{Y}(X)=\max\{\til_{A}(X): A \in Y\}$.
It therefore suffices to show that
$L_{A}(X)=\til_{A}(X)$ for a simplex $A$. Now, $L_A(X) \leq d$ iff $(X,A)$ satisfies $P_d(|A|,0)$, while
$\til_A(X) \leq d$ iff $(X,A)$ satisfies $P_d(0,|A|)$. Finally, $P_d(|A|,0)$ and $P_d(0,|A|)$ are equivalent by
Proposition \ref{p:ppk}.
{\enp}

\section{Empty Intersections and Non-Vanishing Homology}
\label{s:interm}

For a poset $P$ and an element $x \in P$ , let $P_{>x}=\{y \in P:y>x\}$
and $P_{\geq x}=\{y \in P: y\geq x\}$. Let $\Delta(P)$ denote the order complex of $P$, i.e.
the simplicial complex on the vertex set $P$ whose simplices are the chains $x_0<\cdots<x_k$.
Let $M$ be a matroid with rank function $\rho_M$ on the ground set $V$.
Let $\ck(M)$ denote the poset of all flats
$K \neq V$ of $M$ ordered by inclusion, and let $\ck_0(M)=\{K \in \ck(M): \rho_M(K) >0\}$.
It is classically known (see e.g. \cite{Bjorner92}) that $\thh_j\left(\Delta\left(\ck_0(M)\right)\right)=0$
for $j \neq \rho_M(V)-2$. Let $K \in \ck(M)$ and let $B_K$ be an arbitrary basis of $K$. The
contraction of $K$ from $M$ is the matroid on $V \setminus K$ defined by
$M/K=\{A \subset V\setminus K: B_K \cup A \in M\}$ (see e.g. \cite{Oxley11}).
The matroid $M/K$ satisfies $\rho_{M/K}(V \setminus K)= \rho_M(V)-\rho_M(K)$ and
$\ck_0(M/K) \cong \ck(M)_{>K}$.
Let $\{Y_K: Y\in \ck(M)\}$ be a family of simplicial complexes such that
 $Y_K \cap Y_{K'}=Y_{K \cap K'}$ for all $K, K' \in \ck(M)$.
Let $Y=\bigcup_{K \in \ck(M)}Y_K$. For $y \in Y$ let $K_y=\bigcap\{K \in \ck(M): y \in Y_K\} \in \ck(M)$.
The proof of the following result is an application of the method of simplicial resolutions (see e.g. Vassiliev's paper \cite{Vassiliev01}).
\begin{proposition}
\label{p:specs1}
There exists a first quadrant spectral sequence $\{E_{p,q}^r\}$ converging to
$H_*(Y)$ whose $E^1$ term satisfies
\begin{equation}
\label{e:aldu0}
E_{p,q}^1\cong \bigoplus_{\substack{K \in \ck(M) \\ \rho_M(K)=\rho_M(V)-p-1}} H_q(Y_K) \otimes
\thh_{p-1}\big(\Delta\left(\ck_0(M/K)\right)\big).
\end{equation}
\end{proposition}
\noindent
{\bf Proof.}
Let $\rho_M(V)=m$.
For $0 \leq p \leq m-1$ let
$$
F_p=\bigcup_{\substack{K \in \ck(M) \\ \rho_M(K) \geq m-p-1}} Y_K \times
\Delta\big(\ck(M)_{\geq K}\big).
$$
Let $\varphi: F_{m-1}\rightarrow  Y$ denote the projection on the first coordinate. For $y \in Y$, the fiber
$\varphi^{-1}(y)=\{y\} \times \Delta\left(\ck(M)_{\geq K_y}\right)$
is contractible.  Hence, by the Vietoris-Begle theorem (see e.g. p. 344 in \cite{Spanier81}), $H_*(F_{m-1}) \cong H_*(Y)$.
The filtration $F_0 \subset \cdots \subset F_{m-1}$ thus gives rise to a spectral sequence $\{E_{p,q}^r\}$
that converges to $H_*(Y)$.
Let
$$
G_p= \bigcup_{\substack{K \in \ck(M)  \\ \rho_M(K) = m-p-1}} Y_K\times
\Delta\big(\ck(M)_{\geq K}\big).
$$
Then $F_p=G_p \cup F_{p-1}$ and
$$
G_p \cap F_{p-1}=\bigcup_{\substack{K \in \ck(M)  \\ \rho_M(K) = m-p-1}} Y_K\times
\Delta\left(\ck(M)_{> K}\right).
$$
Additionally, if $K \neq K' \in \ck(M)$ satisfy $\rho_M(K) = \rho_M(K')=m-p-1$, then
\begin{equation}
\label{e:addit}
\Big(Y_K\times
\Delta\left(\ck(M)_{\geq K}\right)\Big)
\cap
\Big(Y_{K'}\times
\Delta\left(\ck(M)_{\geq K'}\right)\Big) \subset Y_K\times
\Delta\left(\ck(M)_{> K}\right).
\end{equation}
Using excision, Eq. (\ref{e:addit}), and the K\"{u}nneth formula,  it follows that
\begin{equation}
\label{e:spek1}
\begin{split}
E_{p,q}^1&=H_{p+q}(F_p,F_{p-1})\\
&\cong H_{p+q}(G_p,G_p \cap F_{p-1}) \\
&=H_{p+q}\Big( \bigcup_{\rho_M(K)=m-p-1} Y_K\times
\Delta\left(\ck(M)_{\geq K}\right),  \bigcup_{\rho_M(K)= m-p-1} Y_K \times
\Delta\left(\ck(M)_{> K}\right) \Big) \\
&\cong \bigoplus_{\rho_M(K)=m-p-1} H_{p+q}\Big(Y_K \times \Delta\left(\ck(M)_{\geq K}\right),
Y_K \times \Delta\left(\ck(M)_{> K}\right) \Big)\\
&\cong \bigoplus_{\rho_M(K)=m-p-1} \bigoplus_{i+j=p+q} H_i(Y_K) \otimes
H_j\Big(\Delta\left(\ck(M)_{\geq K}\right), \Delta\left(\ck(M)_{> K}\right) \Big) \\
&\cong \bigoplus_{\rho_M(K)=m-p-1} \bigoplus_{i+j=p+q} H_i(Y_K) \otimes
\thh_{j-1}\big(\Delta\left(\ck(M)_{> K}\right)\big) \\
&\cong \bigoplus_{\rho_M(K)=m-p-1} \bigoplus_{i+j=p+q} H_i(Y_K) \otimes
\thh_{j-1}\big(\Delta\left(\ck_0(M/K)\right)\big).
\end{split}
\end{equation}
\noindent
As $\thh_{j-1}\big(\Delta\left(\ck_0(M/K)\right)\big)=0$
for $j-1 \neq \rho_{M/K}(V \setminus K)-2=\rho_M(V)-\rho_M(K)-2=p-1$, it follows from
(\ref{e:spek1}) that
\begin{equation*}
\label{e:aldu1}
E_{p,q}^1\cong \bigoplus_{\rho_M(K)=m-p-1} H_q(Y_K) \otimes
\thh_{p-1}\big(\Delta\left(\ck_0(M/K)\right)\big).
\end{equation*}
{\enp}
\noindent
For the proof of the next result, we recall the following well-known consequence of the classical Alexander duality, see e.g. section 6 in \cite{Kalai83} and Theorem 2 in \cite{BBM97}. For a simple direct proof see section 2 in \cite{BCP99}.
\begin{theorem}[Combinatorial Alexander Duality]
\label{t:ad1}
Let $X$ be a simplicial complex on $V$. Then for all $0 \leq q \leq |V|-1$
\begin{equation*}
\label{e:adual}
\thh_{|V|-2-q}(X) \cong \thh^{q-1}(X^{\vee}).
\end{equation*}
\end{theorem}
\ \\ \\
Let $\{Z_K: Z\in \ck(M)\}$ be a family of complexes such that
 $Z_K \cup Z_{K'}=Z_{K \cap K'}$ for all $K, K' \in \ck(M)$.
Proposition \ref{p:specs1} implies the following
\begin{corollary}
\label{c:intz}
Suppose that $\bigcap_{K \in \ck(M)}Z_K=\{\emptyset\}$. Then there exist $0 \leq p \leq \rho_M(V)-1$ and $K \in \ck(M)$ of rank $\rho_M(K)=\rho_M(V)-p-1$, such that $\thh_{p-1}(Z_K) \neq 0$.
\end{corollary}
\noindent
{\bf Proof.} Let $\rho_M(V)=m$. We may assume that all the $Z_K$'s are subcomplexes of the $(N-1)$-dimensional simplex
$\Delta_{N-1}$ where $N>m$.
Let $Y_K=Z_K^{\vee}$ be the Alexander dual of $Z_K$ in $\Delta_{N-1}$.
Then for all $K, K' \in \ck(M)$
$$
Y_K \cap Y_{K'}=Z_K^{\vee} \cap Z_{K'}^{\vee}=\left(Z_K \cup Z_{K'}\right)^{\vee}=
Z_{K \cap K'}^{\vee}=Y_{K \cap K'}.
$$
Moreover,
\begin{equation*}
\label{e:aldu1}
Y= \bigcup_{K \in \ck(M)} Y_K =\bigcup_{K \in \ck(M)} Z_K^{\vee} = \left(\bigcap_{K \in \ck(M)} Z_K\right)^{\vee} = \{\emptyset\}^{\vee}=
\partial \Delta_{N-1}
\cong S^{N-2}.
\end{equation*}
By (\ref{e:aldu0}) there exist $0 \leq p \leq m-1$ and $q \geq 0$ such that
$p+q=N-2$, and a flat $K \in \ck(M)$ of rank $\rho_M(K)=m-p-1$ such that $H_q(Y_K) \neq 0$.
Note that $q=N-2-p>m-2-p \geq 0$. By Alexander duality we obtain
\begin{equation*}
\label{e:alexd}
0 \neq H_q(Y_K) =\thh_q(Y_K)=\thh_q(Z_K^{\vee}) \cong \thh_{N-3-q}(Z_K)=\thh_{p-1}(Z_K).
\end{equation*}
{\enp}

\section{A Relative Topological Colorful Helly Theorem}
\label{s:rtch}
{\bf Proof of Theorem \ref{t:rtch}.} Let $M^*=\{\sigma \subset V: \rho_M(V\setminus \sigma)=\rho_M(V)\}$ be the dual matroid of $M$. The rank function of $M^*$ satisfies
$\rho_{M^*}(A)=|A|-\rho_M(V)+\rho_M(V \setminus A)$. For $K \in \ck(M^*)$, we view
the simplices of $X^{\vee} \setminus X^{\vee}[K]$ as a poset ordered by inclusion, and consider its order complex
$Z_K= \Delta\left(X^{\vee} \setminus X^{\vee}[K]\right)$. Then
$$
Z_K \cup Z_{K'}=\Delta\left(X^{\vee} \setminus X^{\vee}[K]\right) \cup \Delta\left(X^{\vee} \setminus X^{\vee}[K']\right)
=\Delta\left(X^{\vee} \setminus X^{\vee}[K\cap K']\right)= Z_{K \cap K'}.
$$
\noindent
Let
$\sd\left(X^{\vee}[V \setminus K]\right)=\Delta\left(X^{\vee}[V \setminus K]\setminus \{\emptyset\}\right)$
denote the barycentric subdivision of $X^{\vee}[V \setminus K]$.
The inclusion map
$$
\sd\left(X^{\vee}[V \setminus K]\right)
\rightarrow \Delta\left(X^{\vee} \setminus X^{\vee}[K]\right)=Z_K
$$
is a homotopy equivalence. Indeed, the retraction
$\Delta\left(X^{\vee} \setminus X^{\vee}[K]\right) \rightarrow
\sd\left(X^{\vee}[V \setminus K]\right)$ is
given by the simplicial map that sends a vertex $\sigma$ of $\Delta\left(X^{\vee} \setminus X^{\vee}[K]\right)$
to the vertex $\sigma \setminus K$ of $\sd\left(X^{\vee}[V \setminus K]\right)$.
It follows that
there is a homotopy equivalence
\begin{equation}
\label{e:hequiv}
Z_K \simeq \sd\left(X^{\vee}[V \setminus K]\right) \simeq  X^{\vee}[V \setminus K].
\end{equation}
Let $\sigma \in X^{\vee}$. Then $V \setminus \sigma \not\in X$, and hence $V \setminus \sigma  \not\in M$. Therefore $\sigma$ does not contain a basis of $M^*$,
 and thus $\sigma  \subset K$ for some $K \in \ck(M^*)$. Hence $\sigma$ is not a vertex of $Z_{K}$.
It follows that
\begin{equation*}
\label{e:inters}
\bigcap_{K \in \ck(M^*)} Z_K = \{\emptyset\}.
\end{equation*}
By Corollary \ref{c:intz} there exist $0 \leq p \leq \rho_{M^*}(V)-1$ and $K \in \ck(M^*)$ such that
\begin{equation}
\label{e:sof1}
\thh_{p-1}(Z_K) \neq 0
\end{equation}
and
\begin{equation}
\label{e:sof2}
\rho_{M^*}(K)=\rho_{M^*}(V)-p-1.
\end{equation}
As $K \in \ck(M^*)$, it follows that $V \setminus K \not\in M$. The assumption
$Y^{\vee} \subset M$ then implies that $V \setminus K \not\in Y^{\vee}$, hence $K \in Y$. Furthermore, (\ref{e:sof2}) is equivalent to
\begin{equation}
\label{e:matrr}
\rho_M(V \setminus K)=|V|-|K|-p-1.
\end{equation}
Using (\ref{e:sof1}),(\ref{e:hequiv}), Alexander duality and (\ref{e:matrr}), we obtain
\begin{equation*}
\label{e:semifin}
\begin{split}
0 &\neq \thh_{p-1}(Z_K)\cong \thh_{p-1}\left(X^{\vee}[V \setminus K]\right) = \thh_{p-1}\left(\lk(X,K)^{\vee}\right) \\
&\cong \thh_{|V|-|K|-p-2}\big(\lk(X,K)\big)=\thh_{\rho_M(V \setminus K)-1}\big(\lk(X,K)\big).
\end{split}
\end{equation*}
As $K \in Y$, it follows from Proposition \ref{p:lerayef} that $\rho_M(V \setminus K) \leq L_Y(X)$.
Finally, $K \in X$ since $\thh_*\left(\lk(X,K)\right) \neq 0$.
{\enp}

\end{document}